\documentclass[12pt]{amsart}
\usepackage{amsmath,amssymb,amsbsy,amsfonts,latexsym,amsopn,amstext,cite,
                                               amsxtra,euscript,amscd,bm}
\usepackage{url}

\usepackage{mathrsfs}

\usepackage{color}
\usepackage[colorlinks,linkcolor=blue,anchorcolor=blue,citecolor=blue,backref=page]{hyperref}

\usepackage{graphics,epsfig}
\usepackage{graphicx}
\usepackage{float}
\usepackage{epstopdf}
\hypersetup{breaklinks=true}
\usepackage{verbatim}
\usepackage{graphicx}
\usepackage{caption}
\usepackage{subcaption}

\usepackage{soul} 
\usepackage{ulem}

\usepackage{float}
\usepackage{epstopdf}
\hypersetup{breaklinks=true}

\usepackage[np]{numprint}
\npdecimalsign{\ensuremath{.}}

\usepackage{bibentry}

\usepackage[english]{babel}
\usepackage{mathtools}
\usepackage{todonotes}
\usepackage{url}
\usepackage[norefs,nocites]{refcheck}

\usepackage{extarrows} 

\newtheorem{thm}{Theorem}

\newtheorem{lem}[thm]{Lemma}

\newtheorem{cor}[thm]{Corollary}

\newtheorem{rem}[thm] {Remark} 

\numberwithin{equation}{section}
\numberwithin{thm}{section}
\numberwithin{table}{section}

\newcommand{\Z}{{\mathbb Z}} 

\newcommand{\R}{{\mathbb R}}

\newcommand{\F}{{\mathbb F}}

\newcommand{\Fq}{{\mathbb F}_q}

\renewcommand{\d}{{\mathrm{d}}}

\newcommand{\supp}{\operatorname{supp}}

  \def \d{\delta} \def \e{\varepsilon}   \def \l{\lambda}

\newfont{\teneufm}{eufm10}
\newfont{\seveneufm}{eufm7}
\newfont{\fiveeufm}{eufm5}
\newfam\eufmfam
     \textfont\eufmfam=\teneufm
\scriptfont\eufmfam=\seveneufm
     \scriptscriptfont\eufmfam=\fiveeufm
\def\fS{{\mathfrak  S}}

\def\sV{{\mathscr V}}

\def \balpha{\bm{\alpha}}
\def \bbeta{\bm{\beta}}

\def \bphi{\bm{\varphi}}

\def\eqref#1{(\ref{#1})}

\def\vec#1{\mathbf{#1}}

\def\cE{\mathcal E}

\def\le{\leqslant}
\def\leq{\leqslant}
\def\ge{\geqslant}
\def\leq{\leqslant}

\def\cC{{\mathcal C}}

\def\cE{{\mathcal E}}

\def\cI{{\mathcal I}}
\def\cJ{{\mathcal J}}

\def\cN{{\mathcal N}}

\def\cR{{\mathcal R}}

\def\cZ{{\mathcal Z}}

\def\ve{\mathbf e}

\def\e{{\mathbf{\,e}}}

\def\eq{{\mathbf{\,e}}_q}

 \def\0{{\mathbf{0}}}

\def\({\left(}
\def\){\right)}
\def\l|{\left|}
\def\r|{\right|}

\def\rf#1{\left\lceil#1\right\rceil}

\def\mand{\qquad \mbox{and} \qquad}

\newif\ifcomment

\usepackage{graphicx}
\begin{document}

\title[Sali{\'e} sums  and modular square roots of shifted primes]{
Shifted bilinear sums of Sali{\'e} sums and the distribution of modular square roots of shifted primes}

\author[I. E. Shparlinski] {Igor E. Shparlinski}
\address{School of Mathematics and Statistics, University of New South Wales, Sydney NSW 2052, Australia}
\email{igor.shparlinski@unsw.edu.au}

\author[Y. Xiao] {Yixiu Xiao}
\address{School of Mathematical Sciences, Shanghai Jiao Tong University, 800 Dongchuan RD, 200240 Shanghai, China}
\email{asaka1312@sjtu.edu.cn}

\begin{abstract} We establish various upper bounds on Type-I and Type-II shifted 
bilinear sums with Sali{\'e} sums modulo a large prime $q$.  We use these bounds to study, for 
fixed integers $a,b\not \equiv 0 \bmod q$,  the distribution of 
solutions to the congruence $x^2 \equiv ap+b \bmod q$, over primes $p\le P$. This is similar to the recently studied
case of $b = 0$, however  the case $b\not \equiv 0 \bmod q$ exhibits 
some new difficulties. 
\end{abstract}


\maketitle

\tableofcontents

\section{Introduction} 

\subsection{Motivation}
Let $q$ be a large prime number and let $\F_q$ denote the finite field of $q$ elements, 
which we assume to be represented by the elements $\{0, 1, \ldots, q-1\}$.
Furthermore, let 
\[
\e(z) = \exp(2 \pi i z) \mand \eq(z) = \exp(2 \pi iz/q).
\] 

We first recall the definition of {\it Sali{\'e} sums\/}.
\[
S(t  ;q)  = \sum_{z \in \F_q^*} \(\frac{z}{q} \) \eq\(t z+  z^{-1}\), 
\]
where $(\frac{x}{q})$ is the Legendre symbol modulo $q \ge 3$ and $z^{-1}$ in the exponent 
is computed modulo $q$. 

Sali{\'e} sums are very close relatives of {\it Kloosterman sums\/}
\[
K(t  ;q)  = \sum_{z \in \F_q^*} \eq\(t z+  z^{-1}\).
\]

Recently there have been several approaches to estimating various bilinear   Type-I and Type-II sums with Kloosterman sums; 
see~\cite{BagShp, BaSh, BFKMM1, FKM, FKMS, KSWX,KMS1,KMS2, LSZ, MQW, Pasc, Shp, ShpZha, Xi-FKM, XuZh} for such bounds, obtained via various techniques, including algebraic geometric methods as well as methods of Fourier analysis and additive combinatorics.
 We also refer to~\cite{BaRa, BFKMM2, GWZ, KoSh, MQW, Wu1, Wu2} for some arithmetic applications of these bounds, such as asymptotic formulas for moments of Dirichlet $L$-functions.

Sali{\'e} sums also appear in a large variety of arithmetic problems~\cite{DFI1,DFI2, DKSZ, Hom, Ngo, Toth,KSSZ,SSZ1, SSZ2}.

Most of these applications  are based on the following well known  link between Sali{\'e} sums  and quadratic congruences. Namely, assuming that $\gcd(q, t) =1$, we have
\begin{equation}
\label{eq:S eval}
  S(t ;q)   =  \varepsilon_q q^{1/2}      \sum_{\substack{x \in \Z_q \\ x^2=t}} \eq \( 2x \)  , 
\end{equation}
where $ \varepsilon_q  = 1$ if $q \equiv 1 \bmod 4$ and  $ \varepsilon_q  = i$  if $q \equiv 3 \bmod 4$,  
see, for example,~\cite[Lemma~12.4]{IwKow}. 
In particular, if $(t/q) = -1$, then  $S(t ;q) =0$. 

\subsection{Set-up} Let $a$, $b$ and $\lambda$ be integers with  $\gcd(a\lambda ,q)=1$.
Define $\mu$ by $\mu \equiv \lambda/2 \bmod q$, $0 \le \mu < q$.
Then, using some elementary transformations and then the identity~\eqref{eq:S eval}, for each pair of 
integers $m$ and $n$ with $\gcd(amn +b, q) = 1$, we derive 
\begin{equation}
\label{eq:red 2 Salie}
\begin{split}
\sum_{ x^2 \equiv amn + b \bmod{q}}\eq\(\lambda x\) & =
\sum_{x^2 \equiv amn + b \bmod{q}} \eq\(2 \mu x\) \\
& =
\sum_{ (\mu x)^2 \equiv\mu^2\(amn + b\) \bmod{q}}\eq\(2 \mu x\) \\
& =
\sum_{ x^2 \equiv\mu^2\(amn + b\) \bmod{q}}\eq\(2  x\) \\
& = \varepsilon_q^{-1}  q^{-1/2}   S(\mu^2\(amn + b\)  ;q). 
\end{split}
\end{equation}

Hence our bounds on sums of  Sali{\'e} sums have an interpretation of 
exponential sums over solutions to congruences $x^2 \equiv amn + b \bmod q$ 
and we present our result in this setting. 

The most challenging question is obtaining nontrivial bounds on  the following bilinear Type-II sums
\begin{equation}
\label{eq:sum W}
W_{a, b, \lambda}\(\balpha, \bbeta; M,N\) = \sum_{m \le  M} \sum_{n \le  N} \alpha_m \beta_n 
   \sum_{\substack{x^2 \equiv amn + b \bmod{q}}}\eq(\lambda x),
\end{equation}
with   $M,N \ge 1$ and some complex weights $\balpha =  \(\alpha_m\) _{m \le  M}$ and 
$\bbeta=\(\beta_n\)_{n \le  N}$. 

Note that in~\eqref{eq:sum W} and in  all similar sums below,
expressions $\sum_{z \le Z}$ mean that the summation is over all positive integers $z \le Z$.

In the case when $\beta_n=1$ identically, such sums are called Type-I sums, however 
for our applications we need to consider slightly more 
general sums.  Namely, let $\cN =  \(N_m\)_{m \le  M}$ be an arbitrary sequence of 
positive integers with $N_m \le N$
for every positive integer   $m \le M$, where  $M$ and $N$ are some real positive parameters. 

We seek a bound on the sums  
\[
    V_{a, b, \lambda}(\balpha;  M,\cN) = \sum_{m \le M} \sum_{n \le N_m} \alpha_m 
   \sum_{\substack{x^2 \equiv amn + b \bmod{q}}} \eq(\lambda x), 
\]
with  some complex weights $\balpha =  \(\alpha_m\) _{m \le  M}$
and the sequence $\cN =  \(N_m\) _{m \le  M}$.  

Finally, it is also interesting and useful for some applications to have bounds on 
``smooth'' sums without weights:
\[
U_{a, b, \lambda}\(m, N\) = \sum_{n \le  N}  
   \sum_{\substack{x^2 \equiv amn + b \bmod{q}}}\eq(\lambda x). 
\]

 For integers $a$ and $b$ and a real $P \ge 2$ we consider the 
sequence  of elements of $\F_q$ which are square roots of $ap+b$ for primes $p \le P$.
More precisely, we are interested in the multiset
\[
\cR_{a,b}(P) = \{x \in \F_q:~ x^2 \equiv ap+b \bmod q, \ p \le P\},
\]
where throughout this work we use $p$ to denote a prime number.

The case of $1 \le a < q$ and $b = 0$ has recently been studied in the series of papers~\cite{DKSZ, KSSZ, SSZ1, SSZ2}. 
However, here we concentrate on the case of  $1 \le b < q$, which requires a very different approach.
In fact our approach also covers the case of $b = 0$ but the previous results from~\cite{DKSZ, KSSZ, SSZ1, SSZ2}
are stronger in this case. 

In particular, we are interested in the distribution of the elements of $\cR_{a,b}(P)$ and thus in 
obtaining an asymptotic formula for the number $N_{a,b}(H,P)$ of elements of $\cR_{a,b}(P)$ which fall in a
given interval $[0, H) \subseteq [0, q-1]$. Thus, we have
\[
N_{a,b}(H,P)= \sum_{p \le P} \# \{x \in[0, H):~x^2 \equiv ap+b \bmod q\}.
\]
In turn, this question naturally leads us to studying exponential sums
\[
S_{a,b, \lambda}(P)=  \sum_{p \leq P} \sum_{\substack{x^2 \equiv ap+b \bmod{q}}} 
   \eq(\lambda x).
\]

We note that   controlling the size of the set $\cR_{1,0}(P) $ is rather difficult, unless one 
assumes the Generalised Riemann Hypothesis (we refer to~\cite{DKSZ,   SSZ2} for a 
discussion of this issue). On the other hand, by a result of Karatsuba~\cite{Kar}, 
for any fixed $\varepsilon > 0$ there exists some $\delta > 0$ such that, 
uniformly over integers $1 \le a,b <q$,  we have 
\begin{equation}
\label{eq:RaP}
\#\cR_{a,b}(P) =  \pi(P) + O\(P^{1-\delta}\)
\end{equation}
provided that $P \ge q^{1/2 + \varepsilon}$, where, as usual $\pi(P)$ denotes  the number of primes $p \le P$
and  throughout  this paper, the implied 
constants in the equivalent notations 
\[
X = O(Y),  \qquad X \ll Y, \qquad Y \gg X
\]
may depend, where obvious, only on $\varepsilon$ and are absolute otherwise. 
We also refer to work of Kerr~\cite{Kerr} on character sums over shifted primes modulo 
a composite number. 

\section{Main results}

\subsection{Bounds of bilinear sums}
We define
\[
    \|\vec{\alpha}\|_2=\(\sum_{m \leq M} |\alpha_m|^2 \)^{1/2}
\mand
    \|\vec{\beta}\|_{\infty}=\max_{n \leq N} |\beta_n|.
\]

We start with the simplest sums $U_{a, b, \lambda}\(m, N\)$. 

 \begin{thm}\label{thm:Smooth sum} For arbitrary  real $N\ge 1$, we have
\[
U_{a, b, \lambda}\(m, N\)  \ll   q^{1/2} \log q , 
\]
uniformly  over  integers $a$, $b$, $m$ and $\lambda$ with $\gcd(am\lambda ,q)=1$. 
\end{thm}  

Using a variation of the argument in~\cite[Appendix~B.3]{DKSZ} we obtain the 
following estimate, which is a full analogue of~\cite[Equation~(B1)]{DKSZ}.

\begin{thm} 
\label{thm:Bilinear Sum I}
 For arbitrary real numbers $M, N\ge 1$, complex weights 
  $\balpha =  \(\alpha_m\) _{m \le  M}$ and  
a sequence $\cN =  \(N_m\)_{m \le  M}$ of real numbers $N_m \le N$,   we have
\[
 \left|V_{a, b,\lambda}(\balpha;  M,\cN)\right| \le  \sqrt{\|\balpha\|_1 \|\balpha\|_2} M^{1/12} N^{7/12} q^{1/4 + o(1)},
\]
uniformly  over  integers $a$, $b$ and $\lambda$ with $\gcd(a\lambda ,q)=1$, 
provided that  
\begin{equation}
\label{eq:Cond MN}
M \le q, \qquad     MN\leq q^{3/2}, \qquad   M \leq N^2.
\end{equation}
\end{thm}   

We note that for $|\alpha_m| \le M^{o(1)}$ and $M \asymp N$, 
Theorem~\ref{thm:Bilinear Sum I} 
gives  a nontrivial bound on sums $V_{a, b,\lambda}(\balpha;  M,\cN)$ already for $M\ge q^{3/7+\varepsilon}$. 

Next, following the ideas of the proof of~\cite[Theorem~2.2]{KSSZ}, we obtain a bound on the   "smoothed''  modification of 
$V_{a, b,\lambda}(\balpha;  M,\cN)$.  That is,  we define
\[
V_{a, b,\lambda}(\balpha, \bphi;  M,\cN) 
= \sum_{m \le M} \sum_{n \in \Z} \alpha_m \varphi_m(n)
   \sum_{\substack{x^2 \equiv amn + b \bmod{q}}} \eq(\lambda x), 
\]
where $\bphi = \(\varphi_m\)_{m \le M}$ is family of smooth functions $\varphi_m$
with supports $\supp \varphi_m \subseteq [0, N_m]$.

\begin{thm} 
\label{thm:Bilinear Sum I Smooth} Let $r \ge 2$ be a fixed integer.
 For arbitrary real numbers $M, N\ge 1$ with $MN \ll q$, complex weights 
  $\balpha =  \(\alpha_m\) _{m \le  M}$ with $\alpha_m \ll 1$, 
a sequence $\cN =  \(N_m\)_{m \le  M}$ of real numbers $N_m \le N$,   
and a   family $\bphi = \(\varphi_m\)_{m \le M}$ of smooth functions $\varphi_m$
with supports $\supp \varphi_m \subseteq [0, N_m]$ and satisfying 
\[
\|\varphi_m^{(j)}\|_\infty \le C_j  N_m^{-j}, \qquad m \le M, \ j =0,1, \ldots, 
\]  
 with some constants $C_j$ depending only on $j$, we have 
\begin{align*}
 &\left|V_{a,b, \lambda}(\balpha, \bphi ; M,\cN)\right | \\
 & \qquad \qquad \le 
\(M^{3/2-1/2r}  N^{1/2 +1/2r}  + M^{1-1/2r}  N^{1/2r} q^{1/2 - 1/4r} \) q^{o(1)} , 
\end{align*}
uniformly  over  integers $a$, $b$ and $\lambda$ with $\gcd(a\lambda ,q)=1$.
\end{thm}

Note that Theorem~\ref{thm:Bilinear Sum I Smooth} is always trivial on the ``diagonal'' $M = N$ (and in fact for all $M \ge N$) 
but is quite strong when $N$ is much larger than $M$.

\begin{rem} 
Yet another improvement of the bounds in~\cite[Equation~(B1)]{DKSZ} and~\cite[Theorem~2.2]{KSSZ} has been outlined in~\cite[Remark~2.3]{KSWX}. However, 
it is not immediately clear how to adapt  argument of~\cite{KSWX} to bilinear sums with shifts.
However,  we discuss an alternative approach to estimating $V_{a, b,\lambda}(\balpha;  M,\cN)$ in 
Appendix~\ref{app:TwistKloost}.
\end{rem}
 
 \begin{thm}\label{thm:Bilinear Sum II} 
 For arbitrary real numbers $M, N\ge 1$, and complex weights 
  $\balpha =  \(\alpha_m\) _{m \le  M}$ and 
$\bbeta=\(\beta_n\)_{n \le  N}$, we have
\begin{align*}
W_{a, b, \lambda}&\( \balpha, \bbeta; M, N \) \\
& \ll \|\vec{\alpha}\|_{2} \, \|\vec{\beta}\|_{\infty} 
\( M^{1/2}N^{1/2} + M^{1/2}N q^{-1/4} + N q^{1/4} (\log q)^{1/2} \), 
\end{align*}
uniformly  over  integers $a$, $b$ and $\lambda$ with $\gcd(a\lambda ,q)=1$. 
\end{thm}

We note that for $M \le q$, by Theorem~\ref{thm:Bilinear Sum II}  we have
\[
\left|W_{a, b, \lambda}\(\balpha, \bbeta; M,N\)\right|   \le  \|\vec{\alpha}\|_2  \|\vec{\beta}\|_{\infty} \(  M^{1/2}N^{1/2} + N  q^{1/4} (\log q)^{1/2}\).
\]

\begin{rem} 
The bound of Theorem~\ref{thm:Bilinear Sum II} is of the same shape 
as the bound of Fouvry,  Kowalski and Michel~\cite[Theorem~1.17]{FKM} and it is very likely 
that it can be derived 
from this result after verifying that Sali{\'e} sum with shifts satisfies the necessary condition 
of~\cite[Theorem~1.17]{FKM}. Our proof of Theorem~\ref{thm:Bilinear Sum II} is based on a different argument  which we hope may be of use for other applications. 
\end{rem} 

 \begin{rem}\label{rem: M vs N}
Clearly the roles of $M$ and $N$ in the bound of Theorem~\ref{thm:Bilinear Sum II} can 
be switched. 
\end{rem}

We recall that for some applications, including this work, one also needs to restrict the summation to
 truncated hyperbolic domains, that is, for points $(u,v) \in \Z^2$ with $uv \le P$, $u \ge U$, $v \ge V$.
 The restrictions  $u \ge U$, $v \ge V$ are needed as it is necessary to stay away from the cusps of the hyperbola to ensure that both summations are long enough.

 Next we define the following sums over a hyperbolic domain
\[
T_{a, b, \lambda}\(\balpha, \bbeta; P, U,V\) = \sum_{\substack{
mn \le P\\m \ge  U, \, n \ge  V}} \alpha_m \beta_n 
   \sum_{\substack{x^2 \equiv amn + b \bmod{q}}}\eq(\lambda x),
\]
with  $P,U,V \ge 1$ and some complex weights $\balpha =  \(\alpha_m\) _{m \ge  U}$ and 
$\bbeta=\(\beta_n\)_{n \ge  V}$. 

\begin{cor}\label{cor:Bilinear sum Hyperb}
For arbitrary  real numbers $P, U, V\ge 1$, and complex weights 
 $\balpha =  \(\alpha_m\) _{m \ge  U}$ and 
$\bbeta=\(\beta_n\)_{n \ge  V}$, satisfying
\[
|\alpha_m|, |\beta_n|  \le q^{o(1)} , 
\]
we have
\[
\left|T_{a, b, \lambda}\(\balpha, \bbeta; P, U,V\) \right| 
        \le \(P V^{-1/2} + Pq^{-1/4} +P q^{1/4}U^{-1/2}\) (Pq)^{o(1)},
\]  
uniformly  over  integers $a$, $b$ and $\lambda$ with $\gcd(a\lambda ,q)=1$. 
\end{cor}

 \begin{rem}\label{rem: weights}
In applications, one usually encounters sums where there are some smooth weights attached to the 
variable $n$, similar to Theorem~\ref{thm:Bilinear Sum I}, and sometimes also to the product $mn$ instead 
of sharp cut-off as in Corollary~\ref{cor:Bilinear sum Hyperb}. For example, such sums appear in the proof 
of Theorem~\ref{thm:SaP-HB} below. Bounds on such sums can be easily reduced to bounds of sums 
$V_{a, b,\lambda}(\balpha;  M,\cN)$, $W_{a, b, \lambda}\( \balpha, \bbeta; M, N \)$ and 
$T_{a, b, \lambda}\(\balpha, \bbeta; P, U,V\)$ via partial summation, without any loses in their strength.
\end{rem}

\subsection{Bounds of sums over primes and distribution of modular 
square roots of shifted primes}
Using the Vaughan identity~\cite{Vau}, we derive  the following bound on the sums $S_{a,b,\lambda}(P)$.   

\begin{thm}
\label{thm:SaP-Vau}
For arbitrary real   $P$, we have
\[    |S_{a,b,\lambda}(P)|  \le   P^{o(1)}
    \begin{cases}  
        P^{13/18}q^{1/4},  & \text{if } P \in (q^{9/10},\ q^{9/8}],\\
        P^{5/6}q^{1/8}, & \text{if } P \in (q^{9/8},\ q^{5/4}], \\
        P^{2/3}q^{1/3}, & \text{if } P \in (q^{5/4},\ q^{3/2}], \\
        P^{5/6}q^{1/12}, & \text{if } P \in (q^{3/2},\ q^2], \\
        Pq^{-1/4}, & \text{if } P > q^2,
    \end{cases}
\]
uniformly  over  integers $a$ and $b$ with $\gcd(a\lambda ,q)=1$. 
\end{thm}

We have presented  Theorem~\ref{thm:SaP-Vau} in the way we prove it. However 
for smaller values of $P\le q^{3/2}$ we use   the more flexible identity of Heath-Brown~\cite{HB} and derive 
a stronger bound.

\begin{thm}
\label{thm:SaP-HB}
For arbitrary real $P$ with   $q^{3/4} \le P  \le q^{3/2}$, we have
\[    |S_{a,b,\lambda}(P)|  \le P^{7/9+o(1)}   q^{1/6}  \]
uniformly  over  integers $a$ and $b$ with $\gcd(a\lambda ,q)=1$. 
\end{thm}

 \begin{rem}\label{rem: Improvement}
Note that the bound of Theorem~\ref{thm:SaP-HB} is always stronger than the bound of Theorem~\ref{thm:SaP-Vau} 
in the full range $q^{3/4} \le P  \le q^{3/2} $ where it applies.  Furthermore, with some additional work, one can improve 
Theorem~\ref{thm:SaP-Vau} for larger values of $P$ as well. However, since we are mostly interested in small values of $P$ 
we do not pursue this here.
\end{rem}

Quite surprisingly, the bound of Theorem~\ref{thm:SaP-HB}, while it is based on very different argument 
and intermediate bounds,  is of exactly the same form the bound of 
Blomer, Fouvry,  Kowalski, Michel and Mili{\'c}evi{\'c}~\cite[Theorem~1.8]{BFKMM2} on sums of 
normalised Kloosterman sums with prime arguments. 

Combining Theorems~\ref{thm:SaP-Vau} and~\ref{thm:SaP-HB} we obtain 
\begin{equation}
\label{eq:Comb bound}   
|S_{a,b,\lambda}(P)|  \le    P^{o(1)} 
    \begin{cases}  
       P^{7/9}q^{1/6} , & \text{if } P \in (q^{3/4},\ q^{3/2}],\\
        P^{5/6}q^{1/12}, & \text{if } P \in (q^{3/2},\ q^2], \\
        Pq^{-1/4}, & \text{if } P > q^2.
    \end{cases}
\end{equation}

We also note that it is very likely that a general bound Fouvry,  Kowalski and Michel~\cite[Theorem~1.5]{FKM} 
on sums of trace functions with prime arguments 
applies to the sums $S_{a,b,\lambda}(P)$ and implies 
\[
 |S_{a,b,\lambda}(P)|  
 \le  
  P^{o(1)}  \begin{cases}  
       P^{11/12}q^{3/48} , & \text{if } P \in  (q^{3/4},\ q],\\ 
               Pq^{-1/48}, & \text{if } P > q, 
    \end{cases}
\]
which, however, is always weaker than a combination of Theorems~\ref{thm:SaP-Vau} and~\ref{thm:SaP-HB}.

Using the bound~\eqref{eq:Comb bound} together with the classical Erd\H{o}s--Tur\'an inequality, which links the discrepancy and exponential sums (see, for instance,~\cite[Theorem~1.21]{DrTi} or~\cite[Theorem~2.5]{KuNi}), we derive:
\begin{cor}
\label{cor:NaHP}  
For arbitrary real $H$ with $q \ge H \ge 1$ and   $P\ge 2$,  we have  
\[
\left|N_{a,b}(H,P)- \frac{H}{q} \#\cR_{a,b}(P) \right| \le  
  P^{o(1)}  \begin{cases}  
       P^{7/9}q^{1/6} , & \text{if } P \in (q^{3/4},\ q^{3/2}],\\
        P^{5/6}q^{1/12}, & \text{if } P \in (q^{3/2},\ q^2], \\
        Pq^{-1/4}, & \text{if } P > q^2,
    \end{cases}
\] 
uniformly  over  integers $a$ and $b$ with $\gcd(a ,q)=1$. 
\end{cor}   

We also appeal to Remark~\ref{rem: Improvement} for possible improvements of Corollary~\ref{cor:NaHP} for $P\ge q$.

Recalling~\eqref{eq:RaP}, we see that the bounds of  Theorems~\ref{thm:SaP-Vau} and~\ref{thm:SaP-HB} and Corollary~\ref{cor:NaHP}
are nontrivial provided $q^A \ge P  \ge q^{3/4 + \varepsilon}$ for any fixed positive $A$ and $\varepsilon$.

\section{Proof of Theorem~\ref{thm:Smooth sum}} 
\label{sec:smooth sum}

Using the periodicity of $\eq(\cdot)$,  we can split the above sum into several 
complete sums running over all $n  \in \F_q$ and one incomplete sum.  
We note that 
\begin{equation}
\label{eq: Vanish N = q}
U_{a, b, \lambda}\(m, q\) = \sum_{n=1}^q  
   \sum_{\substack{x^2 \equiv amn + b \bmod{q}}}\eq(\lambda x)
=   \sum_{x=1}^q\eq(\lambda x) = 0.
\end{equation}  
Hence,  it is enough to show  that   under the assumption
$N \le q$, 
we have
\begin{equation}
\label{eq:U N<q}
U_{a, b, \lambda}\(m, N\) \le   q^{1/2} \log q. 
\end{equation}

We recall that for any $K \le q$ we have 
\begin{equation}
\label{eq:lin sum}
\left| \sum_{k \le K} \eq(hk) \right| \le 
\begin{cases} 
K & \text{if } h = 0, \\
q/h& \text{if } h =1, \ldots, (q-1)/2,
\end{cases}
\end{equation}
see~\cite[Equation~(8.6)]{IwKow}.

Hence, using the standard completing technique,  which is based on  the 
orthogonality of exponential functions and the bound~\eqref{eq:lin sum},  
see, for example,~\cite[Section~12.2]{IwKow},   we write 
\begin{align*}
U_{a, b, \lambda}&\(m, N\) \\
&= \sum_{h=-(q-1)/2}^{(q-1)/2} \frac{1}{q} 
\sum_{n \le N} \eq(hn)
\sum_{t =1}^q 
\sum_{x^2 \equiv amt + b \bmod{q}} 
\eq\(\lambda x - ht\)\\
&\ll  \sum_{h=-(q-1)/2}^{(q-1)/2} \frac{1}{|h|+1} 
\left|  \sum_{ (t,x) \in \cC_m\(\F_q\) } 
\eq\(\lambda x - ht\)\right|, 
\end{align*}
where  $\cC_m\(\F_q\)$ denotes the set of $\F_q$-rational points $(t,x) \in \F_q^2$ 
on the parabola  $\cC_m:~ x^2 = amt + b$ defined 
over $\F_q$. Since it is trivial to show that for $\lambda \in \F_q^*$ and arbitrary $h \in \F_q$,  the linear form 
$\lambda x - ht$ is not constant along $\cC_m$. We can apply the  bound of  Bombieri~\cite[Theorem~6]{Bomb}
on exponential sums along a curve  to the above  sum over $t$ and $x$, and conclude that 
\[
\sum_{ (t,x) \in \cC_m\(\F_q\) } 
\eq\(\lambda x - ht\) \ll q^{1/2}.
\]
Thus
\[
U_{a, b, \lambda}\(m, N\) 
\ll  q^{1/2} \sum_{h= 0}^{(q-1)/2} \frac{1}{h+1}  \ll  q^{1/2} \log q, 
\]
which gives us~\eqref{eq:U N<q} and  concludes the proof.

 \begin{rem}\label{rem: Bombieri vs Gauss}
Theorem~\ref{thm:Smooth sum} admits a more elementary proof based on classical properties of Gauss sums, 
see~\cite[Section~3.4]{IwKow}. We have chosen to appeal to work of  Bombieri~\cite{Bomb}
in preparation to its use in the proof of Theorem~\ref{thm:Bilinear Sum II} below, where it plays 
a crucial role. 
\end{rem}

\section{Proof of Theorem~\ref{thm:Bilinear Sum I}}

\subsection{Preliminary transformations}

This and some other arguments below   are analogues  of those 
in~\cite[Appendix~B.2]{DKSZ}. We also note that due to the vanishing 
complete sums as in~\eqref{eq: Vanish N = q}, we can assume that $N < q$.

Let
\begin{equation}
\label{eq: f_def}
  f(n)=  \sum_{\substack{x^2 \equiv an + b \bmod{q}}}   \eq(\lambda x).
\end{equation}
So we can write  
\[
    V_{a, b, \lambda}(\balpha;  M,\cN)  = \sum_{m \le M} \sum_{n \le N_m} \alpha_m 
   f(mn).
\]

Choose real parameters $D,E \geq 1$ such that
\[
DE \leq N \mand  2DM <q.
\]
We now apply the  technique of~\cite[Section~4]{FouMich} for each $m$ separately with functions $g_m(t)$, 
which are full analogues of $g(t)$ in~\cite{FouMich},  which correspond to the intervals 
$[1, N_m]$ instead of $[1,N]$. 
Namely,   $g_m(t)$  is a trapezoidal continued function with 
$g_m(t) = 1$ for $t \in [1, N_m]$; $g_m(t) = 0$ for $t \in (-\infty, 0] \cup  [N+1, \infty]$
and a linear function in $[0,1]$ and $[N,N+1]$ (with slopes that  ensure the continuity of $g_m(t)$ for all $t \in \R$). In particular, for the Fourier transform 
\[
\widehat g_m(y) = \int_{-\infty}^\infty g(t) \e(-ty) \mathrm{d} t
\]
by~\cite[Equation~(4.1)]{FouMich}, we have 
\[
\widehat g_m(y)  \ll \min\left\{N_m, |y|^{-1}, |y|^{-2}\right\}\}.
\]
Thus we now have full analogues of~\cite[Equation~(4.3)]{FouMich}
and then of~\cite[Equation~(4.4)]{FouMich}, 
which now takes the form 
\begin{align*}
  DE V_{a, b, \lambda}&(\balpha; M, \cN)\\ & \ll \sum_{d \le D} \sum_{m \le M}
  |\alpha_m|
  \sum_{n \le 2N_m}
  \left| \sum_{E < e \leq 2E} \eta(e)  f\(dm\(d^{-1}n+e\)\)  \right| \log q\\
   & \ll \sum_{d \le D} \sum_{m \le M}
    |\alpha_m|
  \sum_{n \le 2N}
  \left| \sum_{E < e \leq 2E}  \eta(e)  f\(dm\(d^{-1}n+e\)\) \right| \log q, 
\end{align*}
where $\eta(e)$ are complex numbers satisfying  $|\eta(e)| \leq 1$ for $E < e \leq 2E$. 
Note that we have also extended the summation over $n$ to the range $n \le 2N$.

Using the H{\"o}lder inequality and the same standard transformations as, for example, 
in~\cite[Section~2.1]{KMS1}, see also~\cite[Equation~(5.8)]{BFKMM1}, we obtain 
the following full analogue of~\cite[Equation~(B.4)]{DKSZ}: 
\begin{align*}
DE V_{a, b, \lambda}(\balpha; M, \cN) 
&\leq \(\|\vec{\balpha}\|_1 \|\vec{\balpha}\|_2\)^{1/2}  (DN)^{3/4} q^{o(1)} \\
&\qquad  \times \left( \sum_{r \in \Fq} \sum_{s \leq 2DM} \left| \sum_{E < e \leq 2E} \eta(e) f\(s(r + e)\) \right|^4 \right)^{1/4} .
\end{align*}

 Expanding the fourth power, the innermost sum in the second factor becomes 
\[
    \sum_{r \in \mathbb{F}_{q}} \sum_{s \leq 2 D M}\left|\sum_{E < e \leq 2E}  \eta(e)  f\(s(r + e)\) \right|^{4}=\sum_{\ve \in \cE} \eta(\ve) \Sigma(\ve), 
\]
where $\cE$ denotes the set of quadruples $\ve=(e_{1}, e_{2}, e_{3}, e_{4}) \in [E,2E]^4$ , 
the coefficients $\eta(\ve)$ satisfy $|\eta(\ve)| \leq 1$ for all $\ve \in \cE$, and  
\begin{align*}
 \Sigma(\ve)=\sum_{r \in \mathbb{F}_{q}} \sum_{s \leq 2 D M} &f\left(s\left(r+e_{1}\right)\right) f\left(s\left(r+e_{2}\right)\right)\\
        &\qquad \quad \times \overline{f\left(s\left(r+e_{3}\right)\right) f\left(s\left(r+e_{4}\right)\right)}.
    \end{align*}

Let $\cE^\Delta$ be the subset of $\cE$, consisting of  tuples $\ve$ with pairs of equal  entries 
(for instance, such as $e_1 = e_2$ and $e_3 = e_4$ or $e_1 = e_3$ and $e_2 = e_4$). For such tuples $\ve$, the following trivial bound holds:
\begin{equation}
\label{eq:bad e}
\sum_{\ve  \in \cE^\Delta} |\Sigma(\ve)| \ll DE^2 M q.
\end{equation}

For $\ve \notin \cE^\Delta$, we complete the sum over $s$ and derive the bound:
\[
\Sigma(\ve) \ll  \max_{t \in \Fq} |\Sigma(\ve, t)| \log q,
\]
where
\begin{align*}
\Sigma(\ve, t) = \sum_{r, s \in \Fq} \eq(st) 
&f\(s(r + e_1)\) f\(s(r + e_2)\) \\
&\qquad \times \overline{f\(s(r + e_3)\) f\(s(r + e_4)\)}.
\end{align*}

\subsection{Bounding  $\Sigma(\ve)$  for ``non-diagonal'' vectors $\ve$}

 We first fix some $t \in \Fq$.  Since $\ve \notin \cE^\Delta$, there is at least one value 
 among $e_{1}$, $e_{2}$, $e_{3}$, $e_{4}$ which is not repeated among other values. 
 Without loss of generality we can assume that 
\[e_{1} \neq e_{2}, e_{3}, e_{4} .
\]

By the definition~\eqref{eq: f_def} of $f(n)$, we have
\[
\Sigma(\ve, t) = \sum_{r, s \in \Fq} \sum_{\substack{x_j^2=as(r+e_j)+b\\j=1,2,3,4}} \eq\(\lambda(x_1+x_2-x_3-x_4)+st\). 
\]

The case of $s=0$ contributes
\[
\sum_{r \in \Fq} \sum_{\substack{x_j^2= b\\j=1,2,3,4}} \eq\(\lambda(x_1+x_2-x_3-x_4)\) \ll q. 
\]
Thus 
\begin{equation}
\label{eq: Sigma*}
\Sigma(\ve, t) = \Sigma^*(\ve, t) + O(q),
\end{equation}
where 
\begin{align*}
 \Sigma^*(\ve, t)     &= \sum_{r\in \Fq} \sum_{s\in \Fq^*} \sum_{\substack{x_j^2=as(r+e_j)+b\\j=1,2,3,4}} \eq\(\lambda(x_1+x_2-x_3-x_4)+st\)\\
    &= \sum_{s\in \Fq^*}\eq(st) 
    \sum_{r\in \Fq}
    \sum_{ x_1^2 = as(r+e_1)+b} \eq(\lambda x_1)\\
   & \qquad  \qquad  \qquad  \qquad  \qquad   \sum_{\substack{x_j^2=as(r+e_j)+b\\j=2,3,4}} 
    \eq\(\lambda(x_2-x_3-x_4)\).
   \end{align*}

We write $u= as(r+e_1)+b$, which runs through $\Fq$ as $r$ runs through $\Fq$ for $s \neq 0$. And change the variables
$x_1 \rightarrow w$. 
We also define 
\[
f_2=e_2-e_1, \qquad f_3=e_3-e_1,  \qquad f_4=e_4-e_1,
\]
and by our assumption   $\ve$, we have $f_2, f_3, f_4 \in \F_q^*$.

Then we can write $ \Sigma^*(\ve, t)$ as 
\begin{align*}
    \Sigma^*(\ve, t)  &= \sum_{s\in \Fq^*}\eq(st) 
    \sum_{u\in \Fq}
    \sum_{w^2 = u} 
    \sum_{\substack{x_j^2 = w^2+af_js\\ j =2,3,4}} 
    \eq\(\lambda(w+x_2-x_3-x_4)\)\\
    &=\sum_{s\in \Fq^*} \eq(st) 
    \sum_{w\in \Fq} 
   \sum_{\substack{x_j^2 = w^2+af_js\\ j =2,3,4}} 
    \eq\(\lambda(w+x_2-x_3-x_4)\). 
\end{align*}

It is clear that $w=0$ contributes $O(q)$.
\begin{equation}
\label{eq: Sigma**}
 \Sigma^*(\ve, t) =  \Sigma^{**}(\ve, t) + O(q), 
\end{equation}
  where 
\begin{align*}
    \Sigma^{**}(\ve, t) & =  \sum_{s\in \Fq^*}\eq(st) 
    \sum_{w\in \Fq^*} 
   \sum_{\substack{x_j^2 = w^2+af_js\\ j =2,3,4}} 
    \eq\(\lambda(w+x_2-x_3-x_4)\)\\
   & =\sum_{s\in \Fq^*}\eq(st) 
    \sum_{w\in \Fq^*} 
    \sum_{\substack{(x_j/w)^2 = 1+af_js/w^2\\j =2,3,4}} 
    \eq\(\lambda(w+x_2-x_3-x_4)\). 
\end{align*}

We now introduce new variables:
\[
\(x_2/w,  x_3/w, x_4/w\) \mapsto \(x,y,z\), 
\]
with which the above sum becomes
\[      \Sigma^{**}(\ve, t)    =\sum_{s\in \Fq^*}\eq(st) 
    \sum_{w\in \Fq^*} 
    \sum_{\substack{x^2 = 1+ af_2 s/w^2\\ y^2 = 1+af_3s/w^2  \\ z^2 = 1+af_4s/w^2}} 
  \eq\(\lambda w(1+x-y-z)\).
\]
Next we make yet another change of variable $s/w^2 \mapsto u$ and, changing the order of summation, derive 
\[     
 \Sigma^{**}(\ve, t)    =\sum_{u\in \Fq^*} 
    \sum_{\substack{x^2 = 1+a f_2u\\ y^2 = 1+af_3u  \\ z^2 = 1+af_4u}} 
        \sum_{w\in \Fq^*} 
  \eq\(tuw^2 + \lambda w(1+x-y-z)\).
\]
Adding the contribution $O(q)$ from the term with $u = 0$ and then another 
 contribution $O(q)$ from the term with $w = 0$,  we have 
\[     
 \Sigma^{**}(\ve, t)   =\sum_{u\in \Fq}
    \sum_{\substack{x^2 = 1+ af_2u\\ y^2 = 1+af_3u  \\ z^2 = 1+af_4u}} 
        \sum_{w\in \Fq} 
  \eq\(tuw^2 + \lambda w(1+x-y-z)\) + O(q).
\]

This sum is identical to the sum in~\cite[Equation~(B.15)]{DKSZ}, which is estimated as
$O(q)$ in~\cite[Section~B.3]{DKSZ}, see~\cite[Equations~(B.19) and~(B.20)]{DKSZ}. 

Thus $ \Sigma^{**}(\ve, t)  \ll q$ and recalling~\eqref{eq: Sigma*} and~\eqref{eq: Sigma**}
we conclude that 
for $\ve \notin \cE^\Delta$, we have
\begin{equation}
\label{eq: Sigma-Bound}
\Sigma(\ve) \ll  q\log q .
\end{equation}

\subsection{Concluding the argument}
The bounds~\eqref{eq:bad e} and~\eqref{eq: Sigma-Bound} are exactly the same 
as bounds in~\cite[Equation~(B.7)]{DKSZ} and in~\cite[Proposition~B.1]{DKSZ}, 
respectively. Hence the same optimisation procedure for the parameters $D$ and $E$, namely, 
the choice
\[
D = \frac{1}{2} M^{-1/3} N^{2/3} \quad \text{and} \quad E = (MN)^{1/3}, 
\]
implies the desired result. 

\section{Proof of Theorem~\ref{thm:Bilinear Sum I Smooth}}

\subsection{Poisson summation}  \label{sec:Poisson}
Exactly as in~\cite[Section~7]{KSSZ} (since this part of the argument does not depend 
on the specific choice of the functions $f_m$, $m \le M$), we have 
\begin{equation}
\begin{split}
\label{eq:V-Poisson}
V_{a,b, \lambda}(\balpha, \bphi ; M,\cN)
& =
\sum_{m \le M} \alpha_m
  \sum_{n \in \mathbb{Z}} \varphi_m(n) f_m(n)\\
& =  \frac{1}{q^{1/2}}\sum_{m \le  M} \alpha_m 
\sum_{n \in \mathbb{Z}}
\widehat{f}_{m}(n)\,\widehat{\varphi}_m\left(\frac{-n}{q}\right), 
\end{split}
\end{equation}
where 
\[
\widehat{f}(x)
=
\frac{1}{q^{1/2}}
\sum_{n \in \F_q}
f(n)\eq(-x n) \quad \text{and}\quad     \widehat{\varphi}_m(\xi) = \int_{-\infty}^\infty \varphi_m(x) \e(-x\xi) \,\mathrm{d} x.
\]

In~\cite{KSSZ} the sum in~\eqref{eq:V-Poisson} is further reduced to Gauss sums. Here we proceed differently and  express $f_m(n)$ in terms of the Sali{\'e} sums. Namely, we use~\eqref{eq:red 2 Salie}
to write 
\begin{align*}
\widehat{f_m}(n)
&=
\frac{1}{q^{1/2}}
\sum_{y \in \mathbb{F}_q}
f_{m}(y)\, e_q(-n y)\\
&=
\frac{1}{q^{1/2}}
\sum_{y \in \mathbb{F}_q}
\varepsilon_q^{-1}  q^{-1/2}   S\!\left(\mu^2(amy + b)  ;q\right)\, e_q(-n y)\\
&=
\frac{1}{q\varepsilon_q}
\sum_{y \in \mathbb{F}_q}
     \sum_{z \in \mathbb{F}_q^*} \left( \frac{z}{q} \right) \e_q\!\left(\mu^2(amy + b) z + z^{-1}\right)\, e_q(-n y)\\
&=
\frac{1}{\varepsilon_q}
     \sum_{z \in \mathbb{F}_q^*} \left( \frac{z}{q} \right) \e_q\!\left(\mu^2b z + z^{-1}\right)\,
     \frac{1}{q}\sum_{y \in \mathbb{F}_q} e_q\!\left((\mu^2amz - n) y\right)\\
&=
\frac{1}{\varepsilon_q}
      \left( \frac{amn}{q} \right) \e_q\!\left(b(am)^{-1}n + \mu^2am n^{-1}\right),
\end{align*}
where, as before, $|\varepsilon_q|=1$ and    $\mu \equiv \lambda/2 \bmod q$, $0 \le \mu < q$.

Substituting this in~\eqref{eq:V-Poisson} and denoting where $\eta \equiv  \mu^2 \bmod q$, we derive 
\begin{equation}
\label{eq: V Salie}
\begin{split}
&V_{a,b, \lambda}(\balpha, \bphi ; M,\cN)\\
&\qquad =\frac{1}{\varepsilon_q q^{1/2}}\sum_{m \le M} \alpha_m \sum_{n \in \mathbb{Z}}
      \left( \frac{amn}{q} \right) \widehat{\varphi}_m\left(\frac{-n}{q}\right)\\
& \qquad \qquad \qquad \qquad \qquad \times  \e_q\!\left(b(am)^{-1}n + \eta am n^{-1}\right)\\
& \qquad\ll q^{-1/2} \sum_{m \le M}  \\
& \qquad \qquad \qquad \ \left|\sum_{n \in \mathbb{Z}}
      \left( \frac{n}{q} \right) \widehat{\varphi}_m\left(\frac{-n}{q}\right) \e_q\!\left(b(am)^{-1}n + \eta am n^{-1}\right)\right|, 
\end{split}
\end{equation}
where we have also recalled our assumption $\alpha_m \ll 1$.

\subsection{Amplification} Define 
\begin{equation}
\label{eq:U-def}
 U= \frac{q}{MN} +1.
\end{equation}
By  our assumption $ MN \ll q$, we have
\begin{equation}
\label{eq:U small}
U \ll \frac{q}{MN} . 
\end{equation} 

For a fixed $m \le  M$, perform the linear shifts $n \mapsto n + u m$  and average
over $1 \le u \le U$. Together with~\eqref{eq: V Salie}, this gives
\begin{equation}
\label{eq: V  Sigmamn}
  V_{a,b, \lambda}(\balpha, \bphi ; M,\cN) \ll \frac{1}{q^{1/2}U}\sum_{m \le M} \sum_{n \in \mathbb{Z}} |\Sigma(m,n, U)|, 
 \end{equation}
   where for a real $W$ we define
\begin{align*}
\Sigma(m,n; W) = \sum_{u \leq W}
    &  \( \frac{n+mu}{q} \) \widehat{\varphi}_m\(\frac{-(n+mu)}{q}\)\\
     & \quad \times  \e_q\!\left(b(am)^{-1}(n+mu) + \eta am (n+mu)^{-1}\right). 
\end{align*}

Clearly, without loss of generality, we can assume that $N_m \ge N/2$ for all $m \le M$.
Then, repeated integration by parts shows that 
\[
\left| \widehat{\varphi}_m(\xi)\right| 
 \ll N_m\min\!\left\{ 1,  (N_m|\xi|)^{-C}\right\} \ll  N \min \left\{ 1, (N |\xi|)^{-C}\right\}
\]
with any constant 
$C>0$. 

Recalling~\eqref{eq:U-def}, we see from~\eqref{eq:U small}  that for any $1\le u \le U$, we have 
\[
\frac{n+m u}{q}=\frac{n}{q}+O\(\frac{MU}{q}\) =\frac{n}{q}+O\(\frac{1}{N}\).
\]
Hence
\[
\widehat{\varphi}_m\!\left(-\frac{n+m u}{q}\right)
\ll \frac{1}{n^C}
\qquad\text{provided}\quad
|n| \gg \frac{q^{1+\varepsilon}}{N},
\]
for any fixed $\varepsilon>0$. Hence we may restrict the summation~\eqref{eq: V  Sigmamn} at 
$|n| \le q^{1+\varepsilon}/N$ and write 
\begin{equation}
\label{eq: V truncated}
  V_{a,b, \lambda}(\balpha, \bphi ; M,\cN) \ll  \widetilde V_{a,b, \lambda}(\balpha, \bphi ; M,N)  
 + MN^{-10} 
\end{equation}
with   
\[
 \widetilde V_{a,b, \lambda}(\balpha, \bphi ; M, N)  = \frac{1}{q^{1/2}U}\sum_{m \le M} \sum_{|n| \le q^{1+\varepsilon}/N} \left| \Sigma(m,n; U) \right| . 
\]

Next, we observe that, 
\[
\frac{\mathrm d}{\mathrm d \xi}
\widehat{\varphi}_m\(\xi\)
= \int_{-\infty}^\infty \varphi_m\left(x\right) \frac{\partial}{\partial \xi} \e (-x\xi)\,\mathrm{d} x \\
=-  \int_{-\infty}^\infty x\varphi_m\left(x\right) e(-x\xi)\,\mathrm{d} x . 
\]
Integrating by parts we obtain 
\begin{equation}
\label{eq:phim by parts}
\left|\frac{\mathrm d}{\mathrm d \xi} \widehat{\varphi}_m\(\xi\)\right| =  
 \frac{1}{|\xi|} \int_{-\infty}^\infty  \frac{\partial}{\partial x}(x\varphi_m(x)) \e(-x\xi)\,\mathrm{d} x.
\end{equation}
Furthermore,  by our assumption on $\varphi_m$   and their derivatives, we have 
\[
\frac{\mathrm d}{\mathrm d x} \(x\varphi_m(x)\) 
=\varphi_m(x)+x \frac{\mathrm d}{\mathrm d x} \varphi_m(x)
\ll 1+x/N_m 
\]
Hence, we derive from~\eqref{eq:phim by parts} that 
\[
\frac{\mathrm d}{\mathrm d \xi}
\widehat{\varphi}_m\(\xi\)
\ll \frac{1}{|\xi|}\(N_m+ N_m^{-1} \int_0^ {N_m} x  \,\mathrm{d} x\)\\
\ll \frac{1}{|\xi|}N_m 
\]
We now see that by the chain rule,
\begin{equation}
\label{eq:phim n+mu bound}
\frac{\partial}{\partial u}
\widehat{\varphi}_m\left(-\frac{n+m u}{q}\right)
\ll N_m\frac{m}{n+mu}\ll \frac{N_m}{u}\ll \frac{N}{u}.
\end{equation}

Let us fix any $U_0$ with 
\[
\Sigma(m,n; U_0)  =  \max_{1 \leq W \leq U} |\Sigma(m,n; W)|.
\]

Recalling~\eqref{eq:phim n+mu bound}, by  partial summation on
$u$
we arrive at
\begin{equation}
\label{eq: Sigma Sigma}
 \Sigma(m,n; U) \ll  | \widetilde\Sigma(m,n; U_0) | N \log q, 
\end{equation}
where 
\begin{align*}
\widetilde \Sigma(m,n; U_0) =
\sum_{u \leq U_0} &
  \left( \frac{nm^{-1}+u}{q} \right) \\
     & \qquad \times 
     \eq\!\left(A(nm^{-1}+u) + \eta a (nm^{-1}+u)^{-1}\right).
\end{align*}
with
\[
A = ba^{-1}  \mand B = \eta a.
\]
Hence, using~\eqref{eq: Sigma Sigma}, we derive 
\begin{equation}
\label{eq: V Sigma U0}
  \widetilde V_{a,b, \lambda}(\balpha, \bphi ; M,N) 
 \ll \frac{N \log q}{q^{1/2}U} \sum_{m  \le M} \sum_{|n| \le q^{1+\varepsilon}/N} \left| \widetilde \Sigma(m,n; U_0)\right|.
\end{equation}

We now change variables:
$nm^{-1} \rightarrow t \bmod q$ and define
\[
\rho(t) =\#\left\{(m,n):~
m  \le M, \  |n| \le q^{1+\varepsilon}/N, \ n m^{-1} \equiv t   \bmod{q}
\right\}.
\]
Then we can rewrite~\eqref{eq: V Sigma U0} as 
\begin{equation}
\label{eq: V  S0}
  \widetilde V_{a,b, \lambda}(\balpha, \bphi ; M,N) \ll\frac{N \log q}{q^{1/2}U}  S(U_0)
\end{equation}
where 
\[
S(U_0) = 
\sum_{t \in \mathbb{F}_q} \rho(t) 
 \left|\sum_{u \leq U_0}
      \left( \frac{t+u}{q} \right) \e_q\!\left(A(t+u) +B (t+u)^{-1}\right)\right|.
\]

We obviously have 
\begin{equation}
\label{eq:I L1}
\sum_{t \in \mathbb{F}_q}
\rho(t) 
\ \ll\
\frac{q^{1+\varepsilon} M}{N}, 
\end{equation}
and as, for example, in~\cite[Equation~(7.6)]{KSSZ},
\begin{equation}
\label{eq:I L2}
        \sum_{t \in \mathbb{F}_q} \rho(t) ^2
\ll \frac{q^{1+2\epsilon}M}{N}.
\end{equation}

By the H\"older inequality,  and then  using~\eqref{eq:I L1} and~\eqref{eq:I L2}, for any integer $r \ge 1$ we have 
\begin{align*}
S(U_0) ^{2r}
& \ll \(\sum_{t \in \mathbb{F}_q} \rho(t)  \)^{2r-2}
\sum_{t \in \mathbb{F}_q} \rho(t) ^2\\
&\qquad \times \sum_{t \in \mathbb{F}_q} \left|\sum_{u \leq U_0}
      \left( \frac{t+u}{q} \right) \e_q\!\left(A(t+u) +B (t+u)^{-1}\right)\right|^{2r}.
\end{align*}

Inserting this in~\eqref{eq: V  S0} and using~\eqref{eq:I L1} and~\eqref{eq:I L2} 
we obtain
\begin{equation}
 \label{eq: V Holder}
  \widetilde V_{a,b, \lambda}(\balpha, \bphi ; M, N) ^{2r} \ll
q^{r-1+ 2r\varepsilon + o(1)}  M^{2r-1} N U^{-2r}  T_r(U_0), 
\end{equation}
where 
\begin{equation}
 \label{eq: Tr}
T_r(U_0)
 =  \sum_{\nu \in \mathbb{F}_q} \left|\sum_{u \leq U_0}
      \left( \frac{\nu+u}{q} \right) \e_q\!\left(A(\nu+u) +B (\nu+u)^{-1}\right)\right|^{2r}.
\end{equation}

\subsection{Application of the Weil bound} 
Expanding the $2r$-th power  in~\eqref{eq: Tr} and changing the order of summation, yields 
\begin{align*}
T_r(U_0)
      & =\sum_{u_1, \ldots, u_{2r}  \leq U_0} \sum_{t \in \mathbb{F}_q} 
      \prod_{i=1}^{2r}
      \left( \frac{t+u_i}{q} \right)\\
      &\qquad \qquad \qquad \times \eq\!\(\sum_{i=1}^{2r}(-1)^{i} \(A(t+u_i) 
      +B (t+u_i)^{-1}\)\).
\end{align*}

Examining the poles of the rational function in the exponent (as a function of $t$)
we see that it is a constant for only $O\(U_0^{r}\)$ choices of 
$1 \leq u_1, \ldots, u_{2r}  \leq U_0$. For the remaining $O\(U_0^{2r}\)$ 
choices we can apply  the Weil bound (see~\cite[Chapter~2, Theorem~2G]{Schm}) and exactly as 
in  the proof of~\cite[Theorem~2.2]{KSSZ} we obtain
\[
T_r(U_0) \ll U_0^{r} q + U_0^{2r} q^{1/2} \le U^{r} q + U^{2r} q^{1/2}.
\]
Substituting these into~\eqref{eq: V Holder} implies that
\begin{align*}
 \widetilde V_{a,b, \lambda}(\balpha, \bphi ; M,N) ^{2r} &\ll q^{r-1+4r\varepsilon+o(1)}M^{2r-1}N\left(q^{1/2}+\frac{q}{U^r}\right) \\
&\ll q^{r-1/2+4r\varepsilon+o(1)}M^{2r-1}N\left(1+\frac{(MN)^{r}}{q^{r-1/2}}\right).
\end{align*}

Recalling~\eqref{eq: V truncated}, since  $\varepsilon>0$ is  arbitrary,  we conclude the proof.

\section{Estimating Type II sums}
\subsection{Proof of Theorem~\ref{thm:Bilinear Sum II}}

By the Cauchy inequality and the assumption on $\balpha$, we have 
\begin{equation}
\label{eq:W and T}
\left|W_{a, b, \lambda}\(\balpha, \bbeta; M,N\)\right|^2   \le  \|\vec{\alpha}\|_{2}^2\, T , 
\end{equation}
where 
\[
T = \sum_{m \le M} \left| \sum_{n \le N} \beta_n 
              \sum_{\substack{x^2 \equiv amn + b \bmod{q}}}\eq(\lambda x) \right|^2.
\]

Clearly, using the periodicity of $\eq(\cdot)$,  we can split the above sum into several 
complete sums running over all $m \in \F_q$  and one incomplete sum. 
Hence  it is enough to show  that   under the assumption
\begin{equation}
\label{eq:M<q}
M \le q,
\end{equation}
we have
\begin{equation}
\label{eq:T M<q}
T \ll \|\vec{\beta}\|_{\infty}^2 \left( MN   + N^2 q^{1/2} \log q \right),
\end{equation} 
which after substitution in~\eqref{eq:W and T} implies the desired bound.

Squaring out, and  using the assumption on $\bbeta$, we derive 
\[   
  T \le     \|\vec{\beta}\|_{\infty}^2   \sum_{n_1,n_2 \leq N} \left|  \sum_{m \le M}   
          \sum_{\substack{x^2 \equiv amn_1 + b \bmod{q} \\ y^2 \equiv amn_2 + b \bmod{q}}} 
         \eq\(\lambda (x-y)\) \right| .
\]

We first estimate the contribution $T_1$ to $T$  from $O\(N(N/q+1)\)$ pairs $(n_1,n_2)$ with $n_1 \equiv n_2 \bmod q$ or 
with $n_1n_2 \equiv 0 \bmod q$ trivially as 
\begin{equation}
\label{eq:Bound T1}
T_1 \ll \|\vec{\beta}\|_{\infty}^2  MN(N/q+1).
\end{equation}

Next we estimate the contribution 
\[
T_2 = \|\vec{\beta}\|_{\infty}^2  \sum_{\substack{n_1,n_2 \le  N\\n_1n_2(n_1 - n_2) \not \equiv 0 \bmod q}} \left|  \sum_{m \le M}   
          \sum_{\substack{x^2 \equiv amn_1 + b \bmod{q} \\ y^2 \equiv amn_2 + b \bmod{q}}} 
         \eq\(\lambda (x-y)\) \right|
\]
from other pairs. 

Recalling our assumption~\eqref{eq:M<q}, and using the standard completing 
technique, see~\cite[Section~12.2]{IwKow}, we see that 
for 
\[
T_2(n_1,n_2) = \sum_{m \leq M} 
\sum_{\substack{x^2 \equiv amn_1 + b \bmod{q} \\ y^2 \equiv amn_2 + b \bmod{q}}}  \eq\(\lambda (x-y)\),  
\] 
as in Section~\ref{sec:smooth sum}, we have 
\begin{align*}
T_2(n_1,n_2)  &= \sum_{h=-(q-1)/2}^{(q-1)/2} \frac{1}{q} 
\sum_{m \le M} \eq(hm)
\sum_{t =1}^q \\
& \qquad \qquad\qquad \times
\sum_{\substack{x^2 \equiv an_1t + b \bmod{q} \\ y^2 \equiv an_2t + b \bmod{q}}} 
\eq\(\lambda (x-y) - ht\) \\
&\leq  2\sum_{h=0}^{(q-1)/2} \frac{1}{q} 
\left| \sum_{m \le M} \eq(hm) \right| \\
& \qquad \qquad\qquad \times
\left| 
\sum_{t =1}^q
\sum_{\substack{x^2 \equiv an_1t + b \bmod{q} \\ y^2 \equiv an_2t + b \bmod{q}}} 
\eq\(\lambda (x-y) - ht\) 
\right|.
\end{align*}

Since $t$ runs through the whole field $\F_q$, for $\gcd(n_1n_2, q) = 1$, eliminating $t$, we have
\begin{equation}
\label{eq:remove t}
\begin{split}
\sum_{t =1}^q &
\sum_{\substack{x^2 \equiv an_1t + b \bmod{q} \\ y^2 \equiv an_2t + b \bmod{q}}} 
\eq\(\lambda (x-y) - ht\) \\
&= \sum_{n_2(x^2-b) \equiv n_1(y^2-b) \bmod{q}}
\eq\(\lambda x - \lambda y - h a^{-1}n_1^{-1} (x^2-b)\).
\end{split}
\end{equation}

We now show that  the bound of  Bombieri~\cite[Theorem~6]{Bomb}
applies to the above  sum over $x$ and $y$, 
provided $n_1n_2(n_1 - n_2) \not \equiv 0 \bmod q$. 

To see this
we first observe that an elementary examination shows  that for $n_1n_2(n_1 - n_2) \not \equiv 0 \bmod q$ 
the polynomial 
\[
n_2(x^2-b) - n_1(y^2-b) \in \F_q[x,y]
\]
is absolutely irreducible in the algebraic closure of $\F_q$. 
Next, we need to show that the polynomial $f(x,y) = \lambda x - \lambda y + h a^{-1}n_1^{-1} (x^2-b)$
is not constant along the curve $n_2(x^2-b) \equiv n_1(y^2-b) \bmod{q}$.

For this we  take any point $(x,y)\in \F_q^2$ on this curve with $y \not \equiv 0 \bmod{q}$ (we 
can appeal to~\cite[Theorem~11.13]{IwKow} to show the existence of such a point when $q$ is large 
enough, but in 
fact for quadrics this can be shown fully elementary). 
Next, we notice that $(x,-y)$ also belongs to this curve, while $f(x,-y) \not \equiv f(x,y) \bmod{q}$
since  $\lambda \not \equiv 0 \bmod{q}$. 

Therefore, for $\gcd(n_1n_2, q) = 1$, the bound $O\(q^{1/2}\)$ of~\cite[Theorem~6]{Bomb} applies to~\eqref{eq:remove t} for every $h = 1, \ldots, q$ and we obtain
\begin{equation}
\label{eq:Bomb}
\begin{split}
\sum_{t =1}^q 
\sum_{\substack{x^2 \equiv an_1t + b \bmod{q} \\ y^2 \equiv an_2t + b \bmod{q}}} 
\eq\(\lambda (x-y) - ht\)  \ll q^{1/2}. 
\end{split}
\end{equation}

 Hence~\eqref{eq:lin sum} and~\eqref{eq:Bomb} imply that for $\gcd(n_1n_2, q) = 1$ we have 
\[
T_2(n_1,n_2) \\
 \ll q^{1/2} \sum_{h=0}^{(q-1)/2} \frac{1}{h+1} \ll  q^{1/2} \log q, 
\]  
and therefore
\begin{equation}
\label{eq:Bound T2}
T_2 \ll  \|\vec{\beta}\|_{\infty}^2  N^2 q^{1/2} \log q.
\end{equation}

From~\eqref{eq:Bound T1} and~\eqref{eq:Bound T2}
\[
T \le T_1 + T_2 
\ll \|\vec{\beta}\|_{\infty}^2 \, MN\left(\frac{N}{q} + 1\right) 
  + \|\vec{\beta}\|_{\infty}^2 \, N^2 q^{1/2} \log q  
\]
which  implies  the bound~\eqref{eq:T M<q}  and
concludes the  proof.

\subsection{Proof of Corollary~\ref{cor:Bilinear sum Hyperb}}
Let $\alpha_m=0$ if $m\le U$ or $m \ge P/V$. Then,
\begin{align*}
& \left|T_{a, b, \lambda}\(\balpha, \bbeta; P, U,V\) \right| 
\\ & \quad =
\sum_{\log U\le j\le\log (P/V)}\sum_{e^j < m \le e^{j+1}}
\sum_{V<n\le  P/m}\alpha_m \beta_n \  \sum_{\substack{x^2 \equiv amn + b \bmod{q}}}\eq(\lambda x).
\end{align*}
Using the completing techniques, 
see~\cite[Section~12.2]{IwKow}, as in the proof of~\cite[Lemma~2.5]{BFGS}, 
for each $j$ the inner sums can be reduced to bounding the sums 
$W_{a, b, \lambda}\(\widetilde \balpha, \widetilde \bbeta; e^{j+1}, P/e^j\)$
with slightly modified weights $\widetilde \balpha$ and $\widetilde \bbeta$ but
still satisfying the condition of  Theorem~\ref{thm:Bilinear Sum II}. 

Hence, the bound of  Theorem~\ref{thm:Bilinear Sum II} yields  
\begin{equation}
\label{eq:dyadic}
\begin{split}
\sum_{e^j < m \le e^{j+1}} & \sum_{V<n\le  P/m}\alpha_m \beta_n \  \sum_{\substack{x^2 \equiv amn + b \bmod{q}}}\eq(\lambda x)\\
& \quad \le  \(P^{1/2} e^{j/2} + P q^{-1/4} +P q^{1/4}e^{-j/2} \)  (Pq)^{o(1)}. 
\end{split}
\end{equation}
Summing over $j$ with $\log U\le j\le\log (P/V)$, we conclude the proof.


\section{Proof of Theorem~\ref{thm:SaP-Vau}}

\subsection{Preliminaries}

As usual, we use $\Lambda(n)$ to denote the  von Mangoldt function
\[\Lambda(n)=
\begin{cases}
\log p &\qquad\text{if $n$ is a power of a prime $p$,} \\
0&\qquad\text{otherwise.}
\end{cases}
\]

It is easy to see that by partial summation it is enough 
to show that 
\begin{equation}
\label{eq: S Lambda}
    \left| \overline  S_{a,b,\lambda}(P)\right| \le   P^{o(1)}
    \begin{cases}  
        P^{13/18}q^{1/4},  & \text{if } P \in (q^{9/10},\ q^{9/8}],\\
        P^{5/6}q^{1/8}, & \text{if } P \in (q^{9/8},\ q^{5/4}], \\
        P^{2/3}q^{1/3}, & \text{if } P \in (q^{5/4},\ q^{3/2}], \\
        P^{5/6}q^{1/12}, & \text{if } P \in (q^{3/2},\ q^2], \\
        Pq^{-1/4}, & \text{if } P > q^2,
    \end{cases}
\end{equation}
for the sum 
\begin{equation}
\label{eq: S Lambda Def}
\overline  S_{a,b,\lambda}(P)=  \sum_{n \leq P} \Lambda(n) f(n), 
\end{equation}
where
\[
  f(n)=  \sum_{\substack{x^2 \equiv an + b \bmod{q}}} 
   \eq(\lambda x).
\]

We first fix some  real numbers $U,V>1$
with $UV\le P$, to be optimised later. 

We need the following result of Vaughan~\cite{Vau}, which is stated
here in the form given in~\cite[Chapter~24]{Dav}.

\begin{lem}
\label{lem:Vau}
For any complex-valued function $f(n)$ and any real numbers $U,V>1$
with $UV\le P$, we have
\[
\overline S_{a,b,\lambda}(P)\ll\Sigma_1+\Sigma_2+\Sigma_3+|\Sigma_4|,
\]
where 
\begin{align*}
\Sigma_1&=\left|\sum_{m\le U}\Lambda(m)f(m)\right|, \\
\Sigma_2&=\sum_{m \le UV}\left|\sum_{n \le P/m}f(mn)\right| \log (UV),\\
\Sigma_3&=\sum_{m\le V}\max_{w\ge 
0}\left|\sum_{w<n \le P/m}f(mn)\right| \log P,\\
\Sigma_4&= \sum_{\substack{mn\le P\\m>U ,\, n>V}}
 \Lambda(m) \sum_{\substack{d\mid n \\d\le V}}\mu(d)f(mn) , 
\end{align*}
where $\mu(d)$ is the M{\"o}bius function. 
\end{lem}

We now need to estimate the corresponding sums 
 $\Sigma_1, \Sigma_2, \Sigma_3$ and $ |\Sigma_4|$ 
 appearing in   Lemma~\ref{lem:Vau}.
 
 For  $\Sigma_1$, using $|f(n)| \le 2$,  we use the trivial bound 
 \begin{equation}
\label{eq: S1}
 \Sigma_1 \le 2\sum_{m\le U}\Lambda(m) \ll U, 
\end{equation}
which follows from the prime number theorem. 

We also note that in the range where the result is trivial  we can always replace $q^{o(1)}$ 
with $P^{o(1)}$.   



\subsection{Bounding $\Sigma_2$} 
We use two different ways to estimate $\Sigma_2$. 

First we assume that 
 \begin{equation}
\label{eq: Assum P UV}
UV\le (2P)^{2/3} \mand 2P \le q^{3/2}.
\end{equation}

We now cover the rest of the summation over $m$ by $O(\log P)$ 
(possibly overlapping)
intervals of the form $M/2 \le m \le M$ with $M \le UV$. 


Next we fix some parameter $W$ with 
 \begin{equation}
\label{eq: W}
W \le q
\end{equation}
 
We first observe that 
\[
M \le W \le q, \qquad   M \le (2P/M)^2 = N^2,   \qquad  MN \le 2P \le q^{3/2},
\]
under the assumptions~\eqref{eq: W} and~\eqref{eq: Assum P UV}. 
Thus  the conditions~\eqref{eq:Cond MN}  hold and we see that by Theorem~\ref{thm:Bilinear Sum I} we have 
\begin{align*}
\sum_{M/2 \le m \le M}\left|\sum_{n \le P/m}f(mn)\right| & \le M^{5/6} (P/M)^{7/12} q^{1/4 + o(1)}\\
& = M^{1/4} P^{7/12} q^{1/4 + o(1)}. 
\end{align*}
Since the above bound is monotonic with respect to $M$, we conclude that for $ M \le  W$ we have 
 \begin{equation}
\label{eq: Medium M}
\sum_{M/2 \le m \le M} \left|\sum_{n \le P/m}f(mn)\right|  \le  P^{7/12+o(1)}  W^{1/4}  q^{1/4}. 
\end{equation}

Finally for $W < M \le UV$, using Corollary~\ref{cor:Bilinear sum Hyperb},  
 we derive 
 \begin{equation}
\label{eq: Large M}
\begin{split}
\sum_{M/2 \le m \le M}& \left|\sum_{n \le P/m}f(mn)\right|\\
 &\le  \(P (P/M) ^{-1/2} + Pq^{-1/4} +P q^{1/4}M^{-1/2}\) (Pq)^{o(1)}\\
&=   \((PM) ^{1/2} + Pq^{-1/4} +P q^{1/4}M^{-1/2}\) (Pq)^{o(1)}\\
&  \le  \((PUV)^{1/2} + Pq^{-1/4} +P q^{1/4}W^{-1/2}\) P^{o(1)}. 
\end{split} 
\end{equation}

Combining the bounds~\eqref{eq: Medium M} and~\eqref{eq: Large M}, we derive  
\[
\Sigma_2 \le \(P^{7/12}  W^{1/4}  q^{1/4} +(PUV)^{1/2} + Pq^{-1/4} +P q^{1/4}W^{-1/2}\) P^{o(1)}. 
\]

Taking 
\[
W =P^{5/9} \le q^{5/6}, 
\]
so that the assumption~\eqref{eq: W} is satisfied, and also noticing that $Pq^{-1/4} \le P^{13/18} q^{1/4}$ 
for $P \le q^{3/2}$, 
 we obtain  
 \begin{equation}
\label{eq: S2}
\Sigma_2\le   \((PUV)^{1/2}   +P^{13/18} q^{1/4}\) P^{o(1)}. 
\end{equation}

We derive another estimate on $\Sigma_2$.
We choose another parameter $ 1 \le \widetilde W \le UV$ and write 
 \begin{equation}
\label{eq: S2 S21 S22}
\Sigma_2  \le  \(\Sigma_{2,1} + \Sigma_{2,2}\)  \log (UV),
\end{equation}
where 
\[
\Sigma_{2,1}=\sum_{m \le  \widetilde W}\left|\,\sum_{n \le P/m}f(mn)\right|
\quad \text{and } \quad
\Sigma_{2,2}=\sum_{ \widetilde W < m \le UV}\left|\,\sum_{n \le P/m}f(mn)\right|. 
\]
Using Theorem~\ref{thm:Smooth sum} for each $m \le  \widetilde W$ we conclude that 
 \begin{equation}
\label{eq: S2-1}
\Sigma_{2,1}  \le   \widetilde W q^{1/2} (Pq)^{o(1)} . 
\end{equation}
On the other hand, splitting the range $ \widetilde W < m \le UV$ into  intervals of the shape $(e^j, e^{j+1}]$ 
as in the proof of  Corollary~\ref{cor:Bilinear sum Hyperb}, 
and then summing the bound~\eqref{eq:dyadic} over $j$ with $\log  \widetilde W\le j\le\log (UV)$, we derive 
 \begin{equation}
\label{eq: S2-2}
\Sigma_{2,2}  \le  \((PUV)^{1/2} + Pq^{-1/4} +P q^{1/4} \widetilde W^{-1/2}\) (Pq)^{o(1)}. 
\end{equation}

Taking 
\[
 \widetilde W = P^{2/3} q^{-1/6}
\]
and substituting~\eqref{eq: S2-1} and~\eqref{eq: S2-2} in~\eqref{eq: S2 S21 S22}, we derive
 \begin{equation}
\label{eq: S2-B}
\Sigma_2\le   \((PUV)^{1/2} + Pq^{-1/4} +P^{2/3} q^{1/3}\) (Pq)^{o(1)}. 
\end{equation}

\subsection{Bounding  $\Sigma_3$}

Using Theorem~\ref{thm:Smooth sum} for each $m \le V$, we conclude that 
 \begin{equation}
\label{eq: S3}
\Sigma_3\le V q^{1/2} P^{o(1)} . 
\end{equation}

\subsection{Bounding $\Sigma_4$}

We now invoke Corollary~\ref{cor:Bilinear sum Hyperb} to deal with $\Sigma_4$
where we take 
\[
\alpha_m =\Lambda(m)= m^{o(1)}
\]
and 
\[
\beta_n =  \sum_{\substack{d\mid n \\d\le 
V}}\mu(d) \ll \sum_{d\mid n} 1 = n^{o(1)}
\]
by the classical divisor bound, see, for example,~\cite[Equation~(1.81)]{IwKow}. 

Hence
 \begin{equation}
\label{eq: S4}
\Sigma_4\le  \(P V^{-1/2} + Pq^{-1/4} +P q^{1/4}U^{-1/2}\)P^{o(1)}. 
\end{equation}





\subsection{Optimisation of $U$ and $V$}
\label{sec:Opt2}
First, we assume that~\eqref{eq: Assum P UV} holds and 
substitute the bounds~\eqref{eq: S1}, \eqref{eq: S2}, \eqref{eq: S3} and~\eqref{eq: S4} 
in Lemma~\ref{lem:Vau}, we obtain
\begin{align*}
        \left|\overline S_{a,b,\lambda}(P)\right| \le 
        &( U+ (PUV)^{1/2}   +P^{13/18} q^{1/4} + Vq^{1/2}+P V^{-1/2}\\
        &\qquad \qquad \qquad \qquad + Pq^{-1/4} +P q^{1/4}U^{-1/2})P^{o(1)}.
    \end{align*}

    Under~\eqref{eq: Assum P UV}, $Pq^{-1/4}$ is dispensable compared to  $ P q^{1/4}U^{-1/2}$. Hence
    \begin{align*}
        \left|\overline S_{a,b,\lambda}(P)\right| \le 
        &( U+ (PUV)^{1/2}   +P^{13/18} q^{1/4} + Vq^{1/2}+P V^{-1/2}\\
        &\qquad \qquad \qquad \qquad +  P q^{1/4}U^{-1/2})P^{o(1)}.
    \end{align*}

We set $UV=(2P)^{2/3}$ to satisfy~\eqref{eq: Assum P UV}.
To balance $P V^{-1/2}$ and $P q^{1/4}U^{-1/2}$, we take 
\[
    U=P^{1/3}q^{1/4} \mand  V=P^{1/3}q^{-1/4},
\]
where  we stipulate  that $P > q^{3/4}$ so that $U \ge V  > 1$. 
Hence
\[
P V^{-1/2}=P q^{1/4}U^{-1/2} = P^{5/6}q^{1/8}.
\]

The terms $U=P^{1/3}q^{1/4}$ and $Vq^{1/2}=P^{1/3}q^{1/4}$ are negligible compared with $P^{13/18} q^{1/4}$. The term $(PUV)^{1/2}=P^{5/6}$ is negligible compared with $P^{5/6}q^{1/8}$.
Therefore,
 \begin{align*}
\left|\overline{S}_{a,b,\lambda}(P)\right|
&\le (P^{13/18}q^{1/4} + P^{5/6}q^{1/8})P^{o(1)} \\
&\le 
\begin{cases} 
P^{13/18+o(1)}q^{1/4}, & P\in(q^{9/10},q^{9/8}]; \\
P^{5/6+o(1)}q^{1/8}, & P\in(q^{9/8},q^{5/4}].
\end{cases}
\end{align*}

For $P> q^{5/4}$, we substitute the bounds~\eqref{eq: S1}, \eqref{eq: S2-B}, \eqref{eq: S3} and~\eqref{eq: S4}
in Lemma~\ref{lem:Vau}, we obtain 
\begin{align*}
        \left|\overline S_{a,b,\lambda}(P)\right| \le 
        &(U+ (PUV)^{1/2}+P^{2/3}q^{1/3} + Vq^{1/2}+P V^{-1/2}\\
        &\qquad \qquad \qquad \qquad + Pq^{-1/4} +P q^{1/4}U^{-1/2})(Pq)^{o(1)}.
    \end{align*}

Comparing the $(PUV)^{1/2}$ (from $\Sigma_2$) and  $P q^{1/4}U^{-1/2}$ (from $\Sigma_4$), we set 
\[
U=P^{1/2}q^{1/4}V^{-1/2}.
\]
(and observe that with  this choice  $U \le P^{2/3}q^{1/3}$). 
Then
\[
(PUV)^{1/2}=P q^{1/4}U^{-1/2}=P^{3/4}V^{1/4}q^{1/8}
\]
 hence
\begin{align*}
&\left|\overline S_{a,b,\lambda}(P)\right| \\
& \qquad \le \(P^{3/4}q^{1/8}V^{1/4} + Vq^{1/2}+P V^{-1/2} +P^{2/3}q^{1/3} + Pq^{-1/4} \)(Pq)^{o(1)}.
 \end{align*}

Comparing the $P^{3/4}q^{1/8}V^{1/4}$  and  $P V^{-1/2}$ (from $\Sigma_4$), we set 
\[
V=P^{1/3}q^{-1/6}.
\]
Then 
\[
P^{3/4}V^{1/4} q^{1/8}=P V^{-1/2}=P^{5/6}q^{1/12},
\]
and $Vq^{1/2} = P^{1/3}q^{1/3}$ is negligible compared to $P^{2/3}q^{1/3}$.
We also have $W =P^{2/3} q^{-1/6}\le UV=P^{2/3}q^{1/6}\leq P$ provided that $P \ge q^{1/2}$.
The upper bound of $|\overline S_{a,b,\lambda}(P)|$ becomes 
\[
    \left|\overline S_{a,b,\lambda}(P)\right| \le \(P^{5/6}q^{1/12}  +P^{2/3}q^{1/3} + Pq^{-1/4} \)(Pq)^{o(1)}.
\]
Noticing that  as before $(Pq)^{o(1)}$ can be replaced with $P^{o(1)}$, we obtain~\eqref{eq: S Lambda}
and conclude the proof.

\section{Proof of Theorem~\ref{thm:SaP-HB}}

\subsection{Preliminaries} 
We again concentrate on the 
sums $\overline S_{a,b,\lambda}(P)$ given by~\eqref{eq: S Lambda Def}. 

 Let 
 \begin{equation}
 \label{eq:cond H1}
 P^{1/2}\ge H\ge P^{\varepsilon}
 \end{equation}
 be some parameter to be chosen later and take  
\begin{equation}
 \label{def:J}
J = \rf{\log P/\log H}  \ge 2.
\end{equation}
All implied constants are allowed to depend on $J$, so we always assume that $H$ exceeds some 
small power of $P$.

Using the Heath-Brown identity~\cite{HB} and a smooth partition of unity, 
see~\cite [Lemma~4.3]{FKM}, exactly as in~\cite[Section~7.1]{DKSZ} or~\cite[Section~8.1]{KSSZ},
Let  $J\ge 1$ be some fixed integer;  in particular, all implied constants are allowed to depend on $J$.  
Then   there exists some $2J$-tuple
\[
\mathbf{V}=(M_1,\ldots , M_J,N_1,\ldots  ,N_J) \in [1/2,2P]^{2J} 
\]
with 
 \begin{equation} \label{eq:size MNQ}
N_1 \geq \ldots \ge N_J, \quad  M_1,\ldots ,M_J \leq P^{1/J},\quad   P \ll  Q \ll  P,
\end{equation}  
where
 \begin{equation} \label{eq:prod Q}
Q =  \prod_{i=1}^J M_i \prod_{j =1}^JN_j,
\end{equation}  
and
\begin{itemize}
\item the arithmetic functions $m_i \mapsto \gamma_i(m_i)$ are bounded and supported in $[M_i/2,2M_i]$;
\item the smooth functions $V_i(x)$, $i = 1, \ldots, J$ and $W(x)$ have support in $[1/2,2]$ and satisfy
\[
V_i^{(j)}(x), W^{(j)}(x) \ll  q^{j \varepsilon}
\]
for all integers $j \geq 0$, where the implied constant may depend on $j$ and $\varepsilon$, 
\end{itemize}  
such that
\[
|\overline S_{a,b,\lambda}(P)| \le \Sigma(\mathbf{V}) P^{o(1)}, 
\]
with 
\begin{align*}
\Sigma(\mathbf{V})& =\sum_{m_1, \ldots, m_J=1}^{\infty} \gamma_1(m_1)\cdots  \gamma_J(m_J)\\
 &\qquad \quad \sum_{n_1,\ldots , n_J=1}^\infty  V_1 \( \frac{n_1}{N_1} \)  \cdots V_J \( \frac{n_J}{N_J} \)   W \( \frac{m_1 \cdots m_J n_1 \cdots  n_J}{P} \)\\
&\qquad \qquad \quad \sum_{x^2\equiv am_1 \cdots m_J n_1 \cdots  n_J+b  \bmod q} \eq(\lambda x).
\end{align*}  

In the following estimates on $\Sigma(\mathbf{V})$, we always appeal to the observation in Remark~\ref{rem: weights}.
Also, the bound  of Theorem~\ref{thm:SaP-Vau} allows us to assume that  
 \begin{equation} \label{eq: P vs Q}
 P \le cq^{3/2}
 \end{equation} 
for some small constant $c> 0$ so that we actually have 
 \begin{equation}  \label{eq: Q vs q}
Q < q^{3/2}.
 \end{equation}

Next, we choose  two arbitrary sets $\cI,\cJ\subseteq \{1,\ldots,J\}$ and write 
 \begin{equation} \label{eq:def MN}
M=\prod_{i\in \cI}M_i\prod_{j\in \cJ}N_j \mand  N=Q/M, 
\end{equation}  
where $Q$ is given by~\eqref{eq:prod Q}, 
which as in~\cite{DKSZ,KSSZ} allows us to treat $\Sigma(\mathbf{V})$ as a bilinear sum. 

Note that by~\eqref{eq:size MNQ}  we have $P/M \le N \ll P/M$ for any $M$ and $N$ as in~\eqref{eq:def MN}.
In particular, using Corollary~\ref{cor:Bilinear sum Hyperb},  we have the following analogue of
the bound in~\cite[Equation~(8.2)]{KSSZ}:
\begin{equation}
\begin{split} 
\label{eq:SigmaV-I} \Sigma(\mathbf{V})
&\le  M^{1/2} 
\(P^{1/2} + M^{-1/2} P q^{-1/4} +  M^{-1} P  q^{1/4}  \) P^{o(1)} \\
&\le 
\(M^{1/2} P^{1/2} +   P q^{-1/4} +  M^{-1/2} P  q^{1/4}  \) P^{o(1)} .
\end{split} 
\end{equation}
However here we also use Theorem~\ref{thm:Bilinear Sum I} which implies, 
under the conditions~\eqref{eq:Cond MN}, that 
\begin{equation}
\label{eq:SigmaV-II} \Sigma(\mathbf{V})
\le  M^{5/6} N^{7/12} q^{1/4}  P^{o(1)} =  N^{-1/4} P^{5/6+o(1)}  q^{1/4}  .
\end{equation}

As in~\cite{DKSZ, KSSZ},  we proceed, forming various set $\cI$ and $\cJ$,  on a case by case basis depending on the size of $N_1$.

\subsection{Small $N_1$} We first consider the case  
\begin{equation}
 \label{eq: small N1}
 N_1  \le H.
\end{equation}

From the definition of $J$ in~\eqref{def:J} and the condition~\eqref{eq:size MNQ}  we see that 
\begin{equation}
 \label{eq: small M}
M_1,\ldots ,M_J \leq  H.
\end{equation}
Here the argument deviates slightly from that in~\cite[Section~8.2]{KSSZ}.
Namely, we choose an additional parameter $T$
to be optimised later. 

We now start multiplying consecutive elements of the sequence 
\[
M_J,\ldots , M_1,N_J,\ldots  ,N_1
\]
until their product $M$ exceeds $T$ for the first time. 
Since by~\eqref{eq: small N1} and~\eqref{eq: small M} each factor is at most $H$, this implies   that we
can  choose two arbitrary sets $\cI, \cJ \subseteq\{1, \ldots, J\}$ such that for 
$M$ and $N$ as in~\eqref{eq:def MN}  we have
\[
T \ll M \ll H T. 
\]
Note that for $T = P^{1/2}$ used in~\cite{DKSZ,KSSZ} one can get tighter 
inequalities with $P^{1/2} \ll M \ll H^{1/2}P^{1/2}$, however the underlying ``switching'' argument 
does not work with an arbitrary $T$.

Hence by~\eqref{eq:SigmaV-I}
\[
\Sigma(\mathbf{V})\le\(H^{1/2} P^{1/2} T^{1/2} +   P q^{-1/4} +   P T^{-1/2} q^{1/4}  \) P^{o(1)} .
\]
 Choosing 
 \begin{equation}
\label{eq:Def T}
T =  H^{-1/2} P^{1/2}   q^{1/4},
\end{equation}
we obtain 
\[
\Sigma(\mathbf{V})\le\(H^{1/4} P^{3/4}  q^{1/8} +   P q^{-1/4}    \) P^{o(1)} .
\]
Furthmore, $P q^{-1/4} \ll H^{1/4} P^{3/4}  q^{1/8}$ provided $P \le q^{3/2}$, hence
\begin{equation}
\label{eq:case1}
\Sigma(\mathbf{V})\le H^{1/4} P^{3/4+o(1)}  q^{1/8} .
\end{equation}

\subsection{Medium $N_1$}
We now  suppose next that
\[
H\le N_1\le T, 
\]
where $T$ is given by~\eqref{eq:Def T}. 

If $N_2 < H$ then we may argue as before. Namely, we  build $M$ starting with $N_1\le T$ and then by multiplying it 
by  other numbers 
\[M_J,\ldots , M_1,N_J,\ldots  ,N_2 \le H,
\]  
in our disposal.  Hence in this case   we obtain  the bound~\eqref{eq:case1} again.

We now suppose 
\[
H\le N_2\le N_1\le T, 
\]
and define  
\[
M=\prod_{i=1}^{J}M_i\prod_{j=3}^{J}N_j  \quad \text{and} \quad N= N_1N_2 ,
\]
Assuming that the conditions~\eqref{eq:Cond MN} hold,  by~\eqref{eq:SigmaV-II} we have 
\begin{equation}
\label{eq:case2}
\Sigma(\mathbf{V})\le H^{-1/2} P^{5/6+o(1)}  q^{1/4}   P^{o(1)} . 
\end{equation}  

In particular, to satisfy $M \le N^2$ and $M\le q$ it is sufficient to request that 
 \begin{equation}
 \label{eq:cond H2}
 H \ge C \max\{P^{1/6}, P^{1/2}q^{-1/2}\}.
 \end{equation}
 with some absolute constant $C>0$.

\subsection{Large $N_1$}
For $N_1 \ge T$, 
    let
\[
M=\prod_{i=1}^{J}M_i\prod_{j=2}^{J}N_j  \quad \text{and} \quad N= N_1 .
\]
Under the conditions~\eqref{eq:Cond MN}, by \eqref{eq:SigmaV-II} we obtain
\[
    \Sigma(\mathbf{V})
\le   N^{-1/4} P^{5/6+o(1)}  q^{1/4}  \le   T^{-1/4} P^{5/6+o(1)}  q^{1/4}  .
\]
Substituting the value of $T$ from~\eqref{eq:Def T} we obtain
 \begin{equation} 
\label{eq: large N_1 tem}
\Sigma(\mathbf{V})
  \le P^{17/24+o(1)}  q^{3/16} H^{1/8}.
\end{equation}
To satisfy~\eqref{eq:Cond MN}, it is sufficient to ensure  that
 \begin{equation}
 \label{eq:cond T large N_1}
T\ge \max\{Q /q, Q^{1/3}\}. 
 \end{equation}
and~\eqref{eq: Q vs q}.

\subsection{Optimisation of $H$}
Putting together all previous bounds~\eqref{eq:case1}, \eqref{eq:case2}  
and~\eqref{eq: large N_1 tem}  results in
\begin{equation*}
\widetilde{S}_{q}(h,P)\le \fS P^{o(1)} ,
\end{equation*}
where 
\[
\fS\le H^{1/4} P^{3/4}  q^{1/8}  
+ H^{-1/2} P^{5/6}  q^{1/4} 
+  P^{17/24}  q^{3/16} H^{1/8}.
\]
As before, to balance $H^{1/4} P^{3/4}  q^{1/8}$ and $H^{-1/2} P^{5/6}  q^{1/4}$, we take 
\[
    H=P^{1/9}q^{1/6},
\]
then
\[
H^{1/4} P^{3/4}  q^{1/8}=H^{-1/2} P^{5/6}  q^{1/4}=P^{7/9}  q^{1/6}.
\]
We also see that the term
\[
 P^{17/24}  q^{3/16} H^{1/8}=P^{13/18}  q^{5/24}
\]
is dominated by $P^{7/9}  q^{1/6}$ provided $P \ge q^{3/4}$. Hence we have 
\[
\fS \ll P^{7/9}  q^{1/6}. 
\]

It remains to verify the assumptions~\eqref{eq:cond H1}, \eqref{eq:cond H2} and~\eqref{eq:cond T large N_1}  
for the above choice of $H$ and $T$.  Indeed, substituting our choice of $H$ in~\eqref{eq:Def T},  we obtain 
\[
T =  \(P^{1/9}  q^{1/6}\)^{-1/2} P^{1/2}   q^{1/4} = P^{4/9}q^{1/6}.
\]
It is now easy to check that~\eqref{eq:cond H1} holds for $P\ge q^{3/4}$  (in fact even for $P \ge c q^{3/7}$), while~\eqref{eq:cond H2} and~\eqref{eq:cond T large N_1} hold 
in the range~\eqref{eq: P vs Q} (in fact even for $P \le c q^{12/7}$), 
which concludes the proof.

\appendix
\section{An alternative approach to Type-I sums}  \label{app:TwistKloost}

Recalling~\eqref{eq:red 2 Salie}, we see that instead of  $V_{a, b, \lambda}(\balpha;  M,\cN)$ it is enough to estimate 
\[\sV_{a, b, \lambda}(\balpha;  M,\cN) = \sum_{m\le M}\alpha_m\sum_{n \le N_m}  S(\mu^2\(amn + b\)  ;q), 
\]
uniformly  over  integers $a$, $b$ and $\lambda$ with $\gcd(a\lambda ,q)=1$.

We again assume that $M \le q$ and deal with $m= q$ trivially.  
Thus we now assume that in fact $M \le q$  and hence  $\gcd(m,q) = 1$ in all summation ranges.

Similarly to our computation in Section~\ref{sec:Poisson}, 
making the change of variable $z \mapsto (a m)^{-1}z \bmod q$  
we have
\begin{align*}
\sV_{a, b, \lambda}&(\balpha;  M,\cN) \\
&   = \sum_{m \le M} \alpha_m \sum_{n \le N_m}  \sum_{z \in \F_q^*} \(\frac{z}{q} \) \eq\((amn + b)z+  z^{-1}\)\\
& =  \sum_{m \le  M} \alpha_m \sum_{n \le N_m}   \sum_{z \in \F_q^*} \(\frac{amz}{q} \) \eq\((amn + b)(a m)^{-1}z+  a mz^{-1}\)\\
& =\sum_{m \le  M} \alpha_m \(\frac{m}{q} \) \sum_{z \in \F_q^*} \(\frac{z}{q} \)   \eq\(b(am)^{-1} z+ a mz^{-1}\)
\sum_{n \le N_m}  \eq\(zn\). 
\end{align*}

It is now convenient to assume that  $\F_q$ is represented by the set 
\[
\{0, \pm 1, \ldots, \pm(q-1)/2\}.
\]

We partition $\F_q^*$ into $2(I+1)$ sets $\cZ_i^{\pm}$, $i =0, \ldots, I$,  with $I = \rf{\log (N/2)}$. We   define $2(I+1)$  sets
 \begin{align*}
\cZ_0^{\pm} &   = \left\{z \in \Z~: ~ 0 < \pm z \le q/N\right \}, \\
\cZ_i^{\pm}  & = \left\{z \in \Z~: ~\min\{q/2, e^{i} q/N\}\ge \pm z > e^{i-1}q/N\right \}.
\end{align*}
Recalling~\eqref{eq:lin sum},  for $z \in\cZ_i^{\pm}$, $i =0, \ldots, I$, we have 
\[
\sum_{n \le N_m } \eq\(n z\) \ll N_m e^{-i} \le N e^{-i}.
\]

Therefore, 
\begin{equation}
\label{eq:V Ti}
\sV_{a, b, \lambda}(\balpha;  M,\cN)  \ll N \sum_{i =0}^I  e^{-i}\(|T_i^+| + |T_i^-|\),
\end{equation}
where $T_i^\pm $  is of the form
\[
T_i^\pm = \sum_{m \le M} \sum_{z \in\cZ_i^{\pm}}   \xi_m \zeta_{z, N}  \eq\(a^{-1}b m^{-1} z+  a mz^{-1}\)
\]
with 
\[ 
\xi_m =  \alpha_m   \(\frac{m}{q} \) , \qquad 1 \le m \le M, 
\]
and 
\[
\zeta_{z} \le 1, \qquad z \in \F_q.
\]

We now recall the discussion in~\cite[Section~9.10]{FKMS}, which shows that~\cite[Theorem~1.3~(2)]{FKMS} applies to sums with  that trace function 
$K(x) =  \eq\(a^{-1}b x+  ax^{-1}\)$. 

Hence, for any integer $r\ge 1$, under the condition
\begin{equation}
\label{eq:range M}
10q^{3/2r} \le M \le q^{1/2 +3/4r}, 
\end{equation}
we have
\begin{align*}
|T_i^\pm| & \le M \(e^{i} q/N\)\( \(e^{i} q/N\)^{-1/2} +
 \(\frac{q^{3/4 + 7/4r}}{M \(e^{i} q/N)\)}\)^{1/2r}\)  q^{o(1)}\\
 & = e^{i/2} M N^{-1/2} q^{1/2+o(1)} \\
 & \qquad  \qquad  \qquad + e^{i(1-1/2r)} M^{1-1/2r} N^{-1+1/2r} 
 q^{1-1/8r + 7/8r^2+o(1)} .
\end{align*} 
Recalling~\eqref{eq:V Ti} we derive that under the condition~\eqref{eq:range M} we have  
\[
|\sV_{a, b, \lambda}(\balpha;  M,\cN) |
\le M N q^{1/2+o(1)} \(N^{-1/2} + \frac{q^{1/2-1/8r + 7/8r^2+o(1)}} {M^{1/2r} N^{1 -1/2r}} \).
\]

 Yet another approach to
 estimating  exponential  sums as $T_i^\pm$, that is,  bilinear sums 
 with $\eq\(u m^{-1} z+  v mz^{-1}\)$ for some 
$u,v \in \F_q^*$ could be via some ideas 
 from~\cite{Bourg, BourgGar}, where   similar sums with $\eq\(u m z+  v m^{-1}z^{-1}\)$
have been extensively investigated. We expect that the methods of these 
works may be useful for the sums $T_i^\pm$ as well.

\section*{Acknowledgement} 

The authors are grateful to Bryce Kerr for very useful discussions and to Philippe Michel for pointing out that the results from~\cite{FKMS} can be applied to 
the sums $T_i^\pm$ in Appendix~\ref{app:TwistKloost}.
 
During the preparation of this work,  I.S. was supported in part by the Australian Research Council Grant DP230100534
and Y.X. by the China Scholarship Council.

\bibliographystyle{plain}

\end{document}